\documentclass[12pt]{article}
\usepackage{amssymb,amsmath,amsfonts,t1enc}
\newcommand{\C}{\mathbb C}
\newcommand{\R}{\mathbb R}
\newcommand{\Z}{\mathbb Z}
\newcommand{\Q}{\mathbb Q}
\newcommand{\F}{\textbf{\textit{F}}}
\newcommand{\K}{\textbf{\textit{K}}}
\newcommand{\wt}{\widetilde}
\begin{document}
\thispagestyle{empty}
\footnotetext{
\par
\noindent
\footnotesize
{\bf 2000 Mathematics Subject Classification.}
{\bf Primary:} 12E99, 12L12.
\par
\noindent
{\bf Key words and phrases:} field endomorphism,
field endomorphism of ${\Q}_p$, field endomorphism of~$\R$,
$p$-adic number that is algebraic over ${\Q}$, real algebraic number.}
\par
\noindent
\centerline{{\large A discrete form of the theorem that each field}}
\centerline{{\large endomorphism of $\R$~(${\Q}_p$) is the identity}}
\vskip 0.5truecm
\centerline{{\large Apoloniusz Tyszka}}
\vskip 0.5truecm
\par
\noindent
{\bf Summary.}
Let $\K$ be a field and $\F$ denote the prime field in $\K$.
Let $\wt{\K}$ denote the set of all
$r \in \K$ for which there exists a finite set $A(r)$ with
$\{r\} \subseteq A(r) \subseteq \K$ such that each
mapping $f:A(r) \to \K$ that satisfies:
~if $1 \in A(r)$ then $f(1)=1$,
~if $a,b \in A(r)$ and $a+b \in A(r)$ then $f(a+b)=f(a)+f(b)$,
~if $a,b \in A(r)$ and $a \cdot b \in A(r)$ then $f(a \cdot b)=f(a) \cdot f(b)$,
~satisfies also $f(r)=r$. Obviously,
each field endomorphism of $\K$ is the identity on $\wt{\K}$.
We prove: $\wt{\K}$ is a countable subfield of $\K$,
if ${\rm char}(\K) \neq 0$ then $\wt{\K}=\F$, $\wt{\C}=\Q$,
if each element of $\K$ is algebraic over $\F=\Q$ then
$\wt{\K}=\{x \in \K: x {\rm ~is~fixed~for~all~automorphisms~of~} \K\}$,
$\wt{\R}$ is equal to the field of real algebraic numbers,
$\wt{{\Q}_p}$ is equal to the field $\{x \in {\Q}_p: x {\rm ~is~algebraic~over~} \Q\}$.
\vskip 0.5truecm
\par
Let $\K$ be a field and $\F$ denote the prime field in $\K$.
Let $\wt{\K}$ denote the set of all $r \in \K$ for which
there exists a finite set $A(r)$ with
$\{r\} \subseteq A(r) \subseteq \K$ such that each
mapping $f:A(r) \to \K$ that satisfies:
\par
\noindent
{\bf (1)}~~~if $1 \in A(r)$ then $f(1)=1$,
\par
\noindent
{\bf (2)}~~~if $a,b \in A(r)$ and $a+b \in A(r)$ then $f(a+b)=f(a)+f(b)$,
\par
\noindent
{\bf (3)}~~~if $a,b \in A(r)$ and $a \cdot b \in A(r)$ then $f(a \cdot b)=f(a) \cdot f(b)$,
\par
\noindent
satisfies also $f(r)=r$. In this situation
we say that $A(r)$ is adequate for $r$. Obviously,
if $f:A(r) \to \K$ satisfies condition {\bf (2)} and $0 \in A(r)$,
then $f(0)=0$.
If $A(r)$ is adequate for $r$ and
$A(r) \subseteq B \subseteq \K$,
then $B$ is adequate for $r$.
We have:
\vskip 0.2truecm
\par
\noindent
{\bf (4)}~~~~~~~~~~~~~~~~~~~~$\wt{\K} \subseteq \widehat{~\K~}:=\bigcap_{\textstyle \sigma\in{\rm End}(\K)}^{}\limits \{x \in \K: \sigma(x)=x \} \subseteq \K,$
\vskip 0.2truecm
\par
\noindent
$\widehat{~\K~}$ is a field.
Let ${\wt{\K}}_n$ ($n=1,2,3,...$) denote the set of all
$r \in \K$ for which there exists $A(r)$ with
$\{r\} \subseteq A(r) \subseteq \K$ such that ${\rm card}(A(r)) \leq n$ and
$A(r)$ is adequate for~$r$.
Obviously,
$${\wt{\K}}_1 \subseteq {\wt{\K}}_2 \subseteq {\wt{\K}}_3 \subseteq ... \subseteq \wt{\K}=\bigcup_{n=1}^\infty {\wt{\K}}_n.$$
{\bf Theorem 1.} $\wt{\K}$ is a subfield of $\K$.
\vskip 0.2truecm
\par
\noindent
{\it Proof.}
We set
$A(0)=\{0\}$ and $A(1)=\{1\}$,
so $0,1 \in \wt{\K}$.
If $r \in \wt{\K}$ then $-r \in \wt{\K}$,
to see it we set
$A(-r)=\{0,-r\} \cup A(r)$.
If $r \in \wt{\K} \setminus \{0\}$
then $r^{-1} \in \wt{\K}$, to see it we set
$A(r^{-1})=\{1,r^{-1}\} \cup A(r)$.
If $r_1,r_2 \in \wt{\K}$ then $r_1+r_2 \in \wt{\K}$,
to see it we set
$A(r_1+r_2)=\{r_1+r_2\} \cup A(r_1) \cup A(r_2)$.
If $r_1,r_2 \in \wt{\K}$ then $r_1 \cdot r_2 \in \wt{\K}$,
to see it we set
$A(r_1 \cdot r_2)=\{r_1 \cdot r_2\} \cup A(r_1) \cup A(r_2)$.
\vskip 0.2truecm
\par
\noindent
{\bf Corollary 1}. If ${\rm char}({\K}) \neq 0$ then $\wt{\K}=\widehat{~\K~}=\F$.
\vskip 0.2truecm
\par
\noindent
{\it Proof.} Let ${\rm char}(\K)=p$.
The Frobenius homomorphism $\K \ni x \to x^p \in \K$
moves all $x \in \K \setminus \F$. It gives $\widehat{~\K~}=\F$,
so by {\bf (4)} and Theorem 1
$\wt{\K}=\widehat{~\K~}=\F$.
\vskip 0.2truecm
\par
\noindent
{\bf Corollary 2.} $\wt{\C}=\widehat{~\C~}=\Q$.
\vskip 0.2truecm
\par
\noindent
{\it Proof.} The author proved (\cite{Tyszka2004}) that
for each $r \in \C \setminus \Q$ there exists a field automorphism
$f:\C \to \C$ such that $f(r) \neq r$.
By this and {\bf (4)} $\wt{\C} \subseteq \widehat{~\C~} \subseteq \Q$,
so by Theorem~1 $\wt{\C}=\widehat{~\C~}=\Q$.
\vskip 0.2truecm
\par
\noindent
{\bf Theorem 2.} For each $n \in \{1,2,3,...\}$ 
${\rm card}(\wt{\K}_n) \leq (n+1)^{n^2+n+1}$, $\wt{\K}$ is countable.
\vskip 0.2truecm
\par
\noindent
{\it Proof.} 
If ${\rm card}(\K)<n$ then
${\rm card}(\widetilde{\K}_n) \leq {\rm card}(\K)<n<(n+1)^{n^2+n+1}$.
In the rest of the proof we assume that ${\rm card}(\K) \geq n$.
Let $r \in {\wt{\K}}_n$ and some
$A(r)=\{r=x_1,...,x_n\}$ is adequate for $r$.
Let also $x_i \neq x_j$ if $i \neq j$.
We choose all formulae
$x_i=1$ ($1 \leq i \leq n$), $x_i+x_j=x_k$, $x_i \cdot x_j=x_k$
($1 \leq i \leq j \leq n$, $1 \leq k \leq n$) that are satisfied
in $A(r)$. Joining these formulae with conjunctions we get
some formula~$\Phi$. Let $V$  denote the set of variables
in $\Phi$, $x_1 \in V$ since otherwise for any $s \in \K \setminus \{r\}$
the mapping $f={\rm id}(A(r)\setminus\{r\}) \cup \{(r,s)\}$
satisfies conditions {\bf (1)}-{\bf (3)} and $f(r) \neq r$.
The formula
$\underbrace{...~\exists x_i~...}_{\textstyle {x_i \in V,~i \neq 1}} \Phi$~~~
is satisfied in~$\K$ if and only if $x_1=r$.
There are $n+1$ possibilities:
$$
1=x_1,~~...,~~1=x_n,~~1 \not\in \{x_1,...,x_n\}.
$$
For each $(i,j) \in \{(i,j): 1 \leq i \leq j \leq n\}$ there are $n+1$
possibilities:
$$
x_i+x_j=x_1,~~...,~~x_i+x_j=x_n,~~x_i+x_j \not\in \{x_1,...,x_n\}.
$$
For each $(i,j) \in \{(i,j): 1 \leq i \leq j \leq n\}$ there are $n+1$
possibilities:
$$
x_i \cdot x_j=x_1,~~...,~~x_i \cdot x_j=x_n,~~x_i \cdot x_j \not\in \{x_1,...,x_n\}.
$$
Since ${\rm card}(\{(i,j): 1 \leq i \leq j \leq n\})=\frac{n^2+n}{2}$ the
number of possible formulae $\Phi$ does not exceed
$(n+1) \cdot (n+1)^\frac{n^2+n}{2} \cdot (n+1)^\frac{n^2+n}{2}=
(n+1)^{n^2+n+1}$.
Thus
${\rm card}({\wt{\K}}_n) \leq (n+1)^{n^2+n+1}$, so
$\wt{\K}=\bigcup_{n=1}^{\infty}\limits {\wt{\K}}_n$ is countable.
\vskip 0.2truecm
\par
\noindent
{\bf Remark 1.} For any field $\K$ the field $\wt{\K}$ is equal to the subfield
of all $x \in \K$ for which $\{x\}$ is existentially $\emptyset$-definable in $\K$.
This gives an alternative proof of Theorems 6 and~7.
\vskip 0.2truecm
\par
Let a field $\K$ extends $\Q$ and each element
of $\K$ is algebraic over $\Q$. R.~M.~Robinson proved (\cite{Robinson}):
if $r \in \K$ is fixed for all automorphisms of $\K$, then there
exist $U(y),V(y) \in {\Q}[y]$ such that $\{r\}$ is definable in $\K$ by
the formula $\exists y~(U(y)=0 \wedge x=V(y))$. Robinson's theorem implies
the next theorem.
\vskip 0.2truecm
\par
\noindent
{\bf Theorem 3.} If a field $\K$ extends $\Q$ and each element
of $\K$ is algebraic over $\Q$, then
$\wt{\K}=\{x \in \K: x {\rm ~is~fixed~for~all~automorphisms~of~} \K\}$.
\vskip 0.2truecm
\par
We use below "bar" to denote the algebraic closure of a field.
\vskip 0.2truecm
\par
\noindent
{\bf Theorem 4} (\cite{Tyszka2005}). If $\K$ is a field and some subfield of $\K$
is algebraically closed, then $\wt{\K}$ is the prime field in $\K$.
\vskip 0.2truecm
\par
\noindent
{\bf Theorem 5.} If a field $\K$ extends $\Q$ and $r \in \wt{\K}$, then
$\{r\}$ is definable in $\K$ by a formula
of the form $\exists x_1 ... \exists x_m T(x,x_1,...,x_m)=0$,
where $m \in \{1,2,3,...\}$ and $T(x,x_1,...,x_m) \in {\Q}[x,x_1,...,x_m]$.
\vskip 0.2truecm
\par
\noindent
{\it Proof.} From the definition of $\wt{\K}$ it follows
that $\{r\}$ is definable in $\K$ by a finite system~{\bf (S)} of polynomial
equations of the form $x_i+x_j-x_k=0$, $x_i \cdot x_j-x_k=0$, $x_i-1=0$,
cf. the proof of Theorem 2.
If~~$\overline{\Q} \subseteq \K$, then by Theorem~4 each element of $\wt{\K}$
is definable in $\K$ by a single equation $x-w=0$, where~$w \in \Q$.
If~~$\overline{\Q} \not\subseteq \K$, then there exists a polynomial
$$
a_{n}x^{n}+a_{n-1}x^{n-1}+...+a_{1}x+a_{0} \in {\Q}[x]~~~~(n \geq 2,~a_n \neq 0)
$$
having no root in $\K$. By this, the polynomial
$$
B(x,y):=a_{n}x^{n}+a_{n-1}x^{n-1}y+...+a_{1}xy^{n-1}+a_{0}y^{n}
$$
satisfies
\begin{equation}
\tag*{{\bf (5)}}
\forall u,v \in \K \left((u=0 \wedge v=0) \Longleftrightarrow B(u,v)=0 \right),
\end{equation}
see \cite{Davis}. Applying {\bf (5)} to {\bf (S)} we obtain
that {\bf (S)} is equivalent to a single polynomial equation.
\vskip 0.2truecm
\par
\centerline{\bf 1. A discrete form of the theorem that each field endomorphism of $\R$ is}
\centerline{\bf the identity}
\vskip 0.2truecm
\par
\noindent
Let ${\R}^{\rm alg}$ denote the field of real algebraic numbers.
\vskip 0.2truecm
\par
\noindent
{\bf Theorem 6.} $\wt{\R}={\R}^{\rm alg}$.
\vskip 0.2truecm
\par
\noindent
{\it Proof.} We prove:
\vskip 0.2truecm
\par
\noindent
{\bf (6)}~~if $r \in {\R}^{\rm alg}$ then $r \in \wt{\R}$.
\vskip 0.2truecm
\par
\noindent
We present three proofs of {\bf (6)}.
\vskip 0.2truecm
\par
\noindent
{\Large \bf (I).} Let $r \in \R$ be an algebraic number
of degree $n$. Thus there exist integers $a_0, a_1, ...,a_n$
satisfying
$$
a_nr^n+...+a_1r+a_0=0
$$
and $a_n \neq 0$.
We choose $\alpha, \beta \in \Q$ such that $\alpha < r < \beta$
and the polynomial
$$
a_nx^n+...+a_1x+a_0
$$
has no roots in $[\alpha, \beta]$ except $r$.
Let $\alpha=\frac{k_1}{k_2}$, $\beta=\frac{l_1}{l_2}$,
where $k_1,l_1 \in \Z$ and $k_2,l_2 \in \Z \setminus \{0\}$.
We put $a={\rm max}\{|a_0|,|a_1|,...,|a_n|,|k_1|,|k_2|,|l_1|,|l_2|\}$.
Then
$$
A(r)=\{\sum_{i=0}^{n} b_ir^i: b_i \in \Z \cap [-a,a] \}
\cup \{ \alpha, r-\alpha, \sqrt{r-\alpha}, \beta, \beta-r, \sqrt{\beta-r} \}
$$
is adequate for $r$. Indeed,
if $f:A(r) \to \R$ satisfies conditions {\bf (1)}-{\bf (3)} then
\vskip 0.2truecm
\par
\centerline{$a_n{f(r)}^n+...+a_1f(r)+a_0=f(a_nr^n+...+a_1r+a_0)=f(0)=0,$}
\vskip 0.2truecm
\par
\noindent
so $f(r)$ is a root of $a_nx^n+...+a_1x+a_0$. Moreover,
\vskip 0.2truecm
\par
\centerline{$f(r)-\alpha=f(r)-f(\alpha)=f(r-\alpha)=
f((\sqrt{r-\alpha})^2)=(f(\sqrt{r-\alpha}))^2 \geq 0$}
\vskip 0.2truecm
\par
\noindent
and
\par
\centerline{$\beta-f(r)=f(\beta)-f(r)=f(\beta-r)=
f((\sqrt{\beta-r})^2)=(f(\sqrt{\beta-r}))^2 \geq 0$.}
\vskip 0.2truecm
\par
\noindent
Therefore, $f(r)=r$.
\vskip 0.2truecm
\par
\noindent
{\Large \bf (II)} (sketch){\bf.} Let $T(x) \in {\Q}[x] \setminus \{0\}$,
$T(r)=0$. We choose $\alpha, \beta \in \Q$ such that
$\alpha < r < \beta$ and $T(x)$ has no roots in $[\alpha,\beta]$
except $r$. Then the polynomial
$$
(1+x^2)^{{\rm deg}(T(x))}
\cdot
T \left(\alpha+\frac{\beta-\alpha}{1+x^2}\right)
\in {\Q}[x]
$$
has exactly two real roots: $x_0$ and $-x_0$.
Thus $x_0^2 \in \wt{\R}$.
By Theorem 1 $\wt{\R}$ is a field, so $\Q \subseteq \wt{\R}$. Therefore,
$r=\alpha+\frac{\beta-\alpha}{1+x_0^2} \in \wt{\R}$.
\vskip 0.2truecm
\par
\noindent
{\Large \bf (III).} The classical Beckman-Quarles theorem states that each
unit-distance preserving mapping from ${\R}^n$ to
${\R}^n$ ($n \geq 2$) is an isometry
(\cite{Beckman-Quarles}-\cite{Benz1994}, \cite{Everling}, \cite{Lester}).
Author's discrete form of this theorem states
that for each $X,Y \in {\R}^n$ ($n \geq 2$) at algebraic
distance there exists a finite set $S_{XY}$ with
$\{X,Y\} \subseteq S_{XY} \subseteq {\R}^n$
such that each unit-distance preserving mapping
$g:S_{XY} \to {\R}^n$ satisfies
$|X-Y|=|g(X)-g(Y)|$ (\cite{Tyszka2000},~\cite{Tyszka2001}).
\vskip 0.2truecm
\par
\noindent
{\sc Case 1:} $r \in {\R}^{\rm alg}$ and $r \geq 0$.
\vskip 0.2truecm
\par
\noindent
The points $X=(0,0) \in {\R}^2$ and $Y=(\sqrt{r},0) \in {\R}^2$
are at algebraic distance $\sqrt{r}$. We consider the finite set
$S_{XY}=\{(x_1,y_1),...,(x_n,y_n)\}$
that exists by the discrete form of the Beckman-Quarles
theorem. We prove that
\begin{eqnarray*}
A(r)=\{0,~1,~r,~\sqrt{r}\} \cup \{x_i:~~1 \leq i \leq n\} \cup \{y_i:~~1 \leq i \leq n \} & \cup & \\
\{x_i-x_j:~~1 \leq i \leq n,~~1 \leq j \leq n\} \cup \{y_i-y_j:~~1 \leq i \leq n,~~1 \leq j \leq n\} & \cup & \\
\{(x_i-x_j)^2:~~1 \leq i \leq n,~~1 \leq j \leq n\} \cup \{(y_i-y_j)^2:~~1 \leq i \leq n,~~1 \leq j \leq n\}
\end{eqnarray*}
is adequate for $r$.
Assume that $f:A(r) \to \R$ satisfies conditions {\bf (1)}-{\bf (3)}.
We show that
$(f,f): S_{XY} \to {\R}^2$ preserves unit distance.
Assume that $|(x_i,y_i)-(x_j,y_j)|=1$,
where $1 \leq i \leq n,~~1 \leq j \leq n$.
Then $(x_i-x_j)^2+(y_i-y_j)^2=1$ and
\begin{eqnarray*}
1=f(1)&=&\\
f((x_i-x_j)^2+(y_i-y_j)^2)&=&\\
f((x_i-x_j)^2)+f((y_i-y_j)^2)&=&\\
(f(x_i-x_j))^2+(f(y_i-y_j))^2&=&\\
(f(x_i)-f(x_j))^2+(f(y_i)-f(y_j))^2&=&
|(f,f)(x_i,y_i)-(f,f)(x_j,y_j)|^2.
\end{eqnarray*}
Therefore, $|(f,f)(x_i,y_i)-(f,f)(x_j,y_j)|=1$.
By the property of $S_{XY}$
$|X-Y|=|(f,f)(X)-(f,f)(Y)|$.
Therefore,
$(0-\sqrt{r})^2+(0-0)^2=$ $|X-Y|^2$ = $|(f,f)(X)-(f,f)(Y)|^2$ = $(f(0)-f(\sqrt{r}))^2+(f(0)-f(0))^2$.
Since $f(0)=0$, we have $r=(f(\sqrt{r}))^2$.
Thus $f(\sqrt{r})=\pm \sqrt{r}$. It implies
$f(r)=f(\sqrt{r} \cdot \sqrt{r})=(f(\sqrt{r}))^2=r$.
\vskip 0.2truecm
\par
\noindent
{\sc Case 2}: $r \in {\R}^{\rm alg}$ and $r<0$.
\vskip 0.2truecm
\par
\noindent
By the proof for case 1 there exists $A(-r)$ that is adequate
for $-r$.
We prove that $A(r)=\{0,~r\} \cup A(-r)$ is adequate for $r$.
Assume that $f:A(r) \to \R$ satisfies conditions {\bf (1)}-{\bf (3)}.
Then $f_{|A(-r)}:A(-r) \to \R$ satisfies conditions {\bf (1)}-{\bf (3)} defined for
$A(-r)$ instead of $A(r)$. Hence $f(-r)=-r$. Since
$0=f(0)=f(r+(-r))=f(r)+f(-r)=f(r)-r$, we conclude that $f(r)=r$.
\vskip 0.2truecm
\par
\noindent
We prove:
\vskip 0.2truecm
\par
\noindent
{\bf (7)}~~if $r \in \wt{\R}$ then $r \in {\R}^{\rm alg}$.
\vskip 0.2truecm
\par
\noindent
Let $r \in \wt{\R}$ and some $A(r)=\{r=x_1,...,x_n\}$
is adequate for $r$.
Let also $x_i \neq x_j$ if $i \neq j$.
We choose all formulae
$x_i=1$ ($1 \leq i \leq n$), $x_i+x_j=x_k$, $x_i \cdot x_j=x_k$ ($1 \leq i \leq j \leq n$,
$1 \leq k \leq n$) that are satisfied
in $A(r)$. Joining these formulae with conjunctions we get
some formula $\Phi$. Let $V$ denote the set of variables
in $\Phi$, $x_1 \in V$ since otherwise for any $s \in \R \setminus \{r\}$
the mapping $f={\rm id}(A(r)\setminus\{r\}) \cup \{(r,s)\}$
satisfies conditions {\bf (1)}-{\bf (3)} and $f(r) \neq r$.
Analogously as in the proof of Theorem 2:
\vskip 0.2truecm
\par
\noindent
{\bf (8)}~~the formula $\underbrace{...~\exists x_i~...}_{\textstyle {x_i \in V,~i \neq 1}} \Phi$~~~
is satisfied in $\R$ if and only if $x_1=r$.
\vskip 0.2truecm
\par
\noindent
The theory of real closed fields is model
complete ([7, {\sc Theorem} 8.6, p.~130]).
The fields $\R$ and ${\R}^{\rm alg}$ are real closed.
Hence ${\rm Th}(\R)={\rm Th}({\R}^{\rm alg})$.
By this, the sentence
$\underbrace{...~\exists x_i~...}_{\textstyle {x_i \in V}} \Phi$ which is
true in $\R$, is also true in ${\R}^{\rm alg}$. Therefore,
for indices~$i$ with $x_i \in V$
there exist $w_i \in {\R}^{\rm alg}$ such that
${\R}^{\rm alg} \models \Phi[x_i \mapsto w_i]$.
Since $\Phi$ is quantifier free, $\R \models \Phi[x_i \mapsto w_i]$.
Thus, by {\bf (8)} $w_1=r$, so $r \in {\R}^{\rm alg}$.
\vskip 0.2truecm
\par
\noindent
{\bf Remark 2.} Similarly to {\bf (7)} the discrete form of
the Beckman-Quarles theorem does not hold for any $X,Y \in {\R}^n$ ($n \geq 2$)
at non-algebraic distance (\cite{Tyszka2000}).
\vskip 0.2truecm
\par
\noindent
{\bf Remark 3.} A well-known result:
\vskip 0.2truecm
\par
\noindent
\centerline{if $f:\R \to \R$ is a field homomorphism, then $f={\rm id}(\R)$
(\cite{Kuczma}-\cite{Lelong-Ferrand})}
\vskip 0.2truecm
\par
\noindent
may be proved geometrically as follows.
If $f: \R \to \R$ is a field homomorphism
then $(f,f): {\R}^2 \to {\R}^2$ preserves unit distance;
we prove it analogously as in {\Large {\bf (III)}}.
By the classical Beckman-Quarles theorem $(f,f)$ is an isometry.
Since the isometry $(f,f)$ has three non-collinear fixed points:
$(0,0)$, $(1,0)$, $(0,1)$, we conclude that
$(f,f)={\rm id}({\R}^2)$ and $f={\rm id}({\R})$.
\vskip 0.2truecm
\par
\noindent
\centerline{\bf 2. A discrete form of the theorem that each field endomorphism of ${\Q}_p$ is}
\centerline{\bf the identity}
\vskip 0.2truecm
~~Let ${\Q}_p$ be the field of $p$-adic numbers,
${|\cdot|}_p$ denote the $p$-adic norm on ${\Q}_p$,
${\Z}_p=\{x \in {\Q}_p: {|x|}_p \leq 1\}$.
Let $v_p: {\Q}_p \to \Z \cup \{\infty\}$ denote the valuation
function written additively:
$v_p(x)=-{\rm log}_{p}({|x|}_p)$ if $x \neq 0$, $v_p(0)=\infty$.
For $n \in \Z$, $a,b \in {\Q}_p$ by $a \equiv b~({\rm mod}~p^n)$ we
understand ${|a-b|}_p \leq p^{-n}$.
It is known (\cite{Lang1984},\cite{Robert},\cite{Wagner})
that each field automorphism of ${\Q}_p$ is the identity.
\vskip 0.2truecm
\par
\noindent
{\bf Lemma 1} (Hensel's lemma, \cite{Koblitz}).
Let $F(x)=c_0+c_1x+...+c_nx^n \in {\Z}_p[x]$. 
Let $F'(x)=c_1+2c_2x+3c_3x^2+...+nc_nx^{n-1}$ be the formal
derivative of $F(x)$. Let $a_0 \in {\Z}_p$ such that
$F(a_0) \equiv 0~({\rm mod}~p)$ and $F'(a_0) \not\equiv 0~({\rm mod}~p)$.
Then there exists a unique $a \in {\Z}_p$ such that $F(a)=0$ and
$a \equiv a_0~({\rm mod}~p)$.
\vskip 0.2truecm
\par
\noindent
{\bf Lemma 2} (\cite{Delon}). For each $x \in {\Q}_p$ ($p \neq 2$)
${|x|}_p \leq 1$ if and only if there exists $y \in {\Q}_p$
such that $1+p{x}^2=y^2$. For each $x \in {\Q}_2$ ${|x|}_2 \leq 1$
if and only if there exists $y \in {\Q}_2$ such that $1+2{x}^3=y^3$.
\vskip 0.2truecm
\par
\noindent
{\it Proof in case $p \neq 2$.} 
If ${|x|}_p \leq 1$ then $v_p(x) \geq 0$
and $x \in \Z_p$. We apply Lemma 1 for $F(y)=y^2-1-px^2$ and $a_0=1$.
This $a_0$ satisfies the assumptions:
$F(a_0)=-px^2 \equiv 0~({\rm mod}~p)$ and
$F'(a_0)=2 \not\equiv 0~({\rm mod}~p)$.
By Lemma 1 there exists $y \in {\Z}_p$ such that $F(y)=0$, so 
$1+px^2=y^2$. If ${|x|}_p>1$ then $v_p(x)<0$. By this
$v_p(1+px^2)=v_p(px^2)=1+2v_p(x)$ is not divisible by $2$, so
$1+px^2$ is not a square.
\vskip 0.2truecm
\par
\noindent
{\it Proof in case $p=2$.} 
If ${|x|}_2 \leq 1$ then $v_2(x) \geq 0$ and $x \in {\Z}_2$.
We apply Lemma 1 for $F(y)=y^3-1-2x^3$ and $a_0=1$. This $a_0$
satisfies the assumptions:
$F(a_0)=-2x^3 \equiv 0~({\rm mod}~2)$ and $F'(a_0)=3 \not\equiv 0~({\rm mod}~2)$.
By Lemma 1 there exists $y \in \Z_2$ such that $F(y)=0$, so $1+2x^3=y^3$.
If ${|x|}_2>1$ then $v_2(x)<0$. By this
$v_2(1+2x^3)=v_2(2x^3)=1+3v_2(x)$ is not divisible by $3$, so
$1+2x^3$ is not a cube.
\vskip 0.2truecm
\par
\noindent
{\bf Lemma 3.} If $c,d \in {\Q}_p$ and $c \neq d$,
then there exist $m \in \Z$ and $u \in \Q$
such that ${|\frac{c-u}{p^{m+1}}|}_p \leq 1$
and ${|\frac{d-u}{p^{m+1}}|}_p>1$.
\vskip 0.2truecm
\par
\noindent
{\it Proof.} Let $c=\sum_{k=s}^{\infty}\limits {c}_{k}{p}^{k}$
and $d=\sum_{k=s}^{\infty}\limits {d}_{k}{p}^k$,
where $s \in \Z$, $c_{k},d_{k} \in \{0,1,...,p-1\}$.
Then $m={\rm min}\{k: c_{k} \neq d_{k}\}$
and $u=\sum_{k=s}^{m}\limits c_{k}{p}^k$ satisfy our conditions.
\vskip 0.2truecm
\par
\noindent
Let ${\Q}_p^{\rm alg}=\{x \in {\Q}_p: x {\rm ~is~algebraic~over~} \Q\}$.
\vskip 0.2truecm
\par
\noindent
{\bf Theorem 7.} $\wt{{\Q}_p}={\Q}_p^{\rm alg}$.
\vskip 0.2truecm
\par
\noindent
{\it Proof.} We prove: if $r \in {\Q}_p^{\rm alg}$ then $r \in \wt{{\Q}_p}$.
\vskip 0.2truecm
\par
\noindent
Let $r \in {\Q}_p^{\rm alg}$.
Since $r \in {\Q}_p$ is algebraic over $\Q$,
it is a zero of a polynomial $p(x)=a_nx^n+...+a_1x+a_0 \in \Z[x]$
with $a_n \neq 0$. Let $R=\{r=r_1,r_2,...,r_k\}$ be the set of all
roots of $p(x)$ in ${\Q}_p$.
For each $j \in \{2,3,...,k\}$ we apply Lemma 3 for $c=r$ and $d=r_j$
and choose $m_j \in \Z$ and $u_j \in \Q$ such that
${|\frac{r-u_j}{p^{m_j+1}}|}_p \leq 1$ and ${|\frac{r_j-u_j}{p^{m_j+1}}|}_p>1$.
Let $u_j=\frac{s_j}{t_j}$, where $s_j \in \Z$ and
$t_j \in \Z \setminus \{0\}$.
In case $p \neq 2$ by Lemma 2 for each $j \in \{2,3,...,k\}$ there exists $y_j \in {\Q}_p$
such that
$$
1+p{\left(\frac{r-u_j}{p^{m_j+1}}\right)}^2=y_j^2.
$$
In case $p=2$ by Lemma 2 for each $j \in \{2,3,...,k\}$ there exists $y_j \in {\Q}_2$
such that
$$
1+2{\left(\frac{r-u_j}{2^{m_j+1}}\right)}^3=y_j^3.
$$
Let $a={\rm max~}\{p, |a_i|, |s_j|, |t_j|, |m_j+1|:~0 \leq i \leq n,~2 \leq j \leq k \}$. The set
$$
A(r)=\{\sum_{i=0}^n b_ir^i: b_i \in \Z \cap [-a,a] \} \cup \{p^w:~w \in \Z \cap [-a,a] \}~\cup
$$ 
\par
\noindent
\centerline{
$\bigcup_{j=2}^k\limits
\{u_j,~r-u_j,~\frac{r-u_j}{p^{m_j+1}},~
\left(\frac{r-u_j}{p^{m_j+1}}\right)^2,~
p\left(\frac{r-u_j}{p^{m_j+1}}\right)^2,~
\left(\frac{r-u_j}{p^{m_j+1}}\right)^3,~
p\left(\frac{r-u_j}{p^{m_j+1}}\right)^3,~y_j,~y_j^2,~y_j^3 \}$}
\vskip 0.2truecm
\par
\noindent
is finite, $r \in A(r)$. We prove that $A(r)$ is adequate for $r$. Assume that
$f:A(r) \to {\Q}_p$ satisfies conditions {\bf (1)}-{\bf (3)}.
Analogously as in {\Large \bf (I)} we conclude
that $f(r)=r_j$ for some $j \in \{1,2,...,k\}$.
Therefore, $f(r)=r$ if $k=1$. Let $k \geq 2$.
Suppose, on the contrary, that
\vskip 0.2truecm
\par
\noindent
{\Large \bf ($\ast$)}~~~~~~~~~~~~~~~~~~~~~~~~~~$f(r)=r_j$ for some $j \in \{2,3,...,k\}$.
\vskip 0.2truecm
\par
\noindent
In case $p \neq 2$ supposition {\Large \bf ($\ast$)} implies:
$$
1+p{\left(\frac{r_j-u_j}{p^{m_j+1}}\right)}^2
=1+p{\left(\frac{f(r)-u_j}{p^{m_j+1}}\right)}^2
=f\left(1+p{\left(\frac{r-u_j}{p^{m_j+1}}\right)}^2\right)
=f(y_j^2)={f(y_j)}^2.
$$
Thus, by Lemma 2 ${|\frac{r_j-u_j}{p^{m_j+1}}|}_p \leq 1$, a contradiction.
In case $p=2$ supposition {\Large \bf ($\ast$)} implies:
$$
1+2{\left(\frac{r_j-u_j}{2^{m_j+1}}\right)}^3
=1+2{\left(\frac{f(r)-u_j}{2^{m_j+1}}\right)}^3
=f\left(1+2{\left(\frac{r-u_j}{2^{m_j+1}}\right)}^3\right)
=f(y_j^3)={f(y_j)}^3.
$$
Thus, by Lemma 2 ${|\frac{r_j-u_j}{2^{m_j+1}}|}_2 \leq 1$, a contradiction.
\vskip 0.2truecm
\par
\noindent
We prove: if $r \in \wt{{\Q}_p}$ then $r \in {\Q}_p^{\rm alg}$.
\vskip 0.2truecm
\par
\noindent
Let $r \in \wt{{\Q}_p}$ and some $A(r)=\{r=x_1,...,x_n\}$ is
adequate for $r$. Let also $x_i \neq x_j$ if $i \neq j$.
Analogously as in the proof of~{\bf (7)} we construct a quantifier free
formula $\Phi$ such that
\vskip 0.2truecm
\par
\noindent
{\bf (9)}~~the formula
$\underbrace{...~\exists x_i~...}_{\textstyle {x_i \in V,~i \neq 1}} \Phi$~~
is satisfied in ${\Q}_p$ if and only if $x_1=r$;
\vskip 0.2truecm
\par
\noindent
as previously, $V$ denote the set of variables in $\Phi$ and $x_1 \in V$.
${\rm Th}({\Q}_p)={\rm Th}({\Q}_p^{\rm alg})$, it follows from the first sentence
on page 134 in \cite{MacIntyre1986}, see also [14, Theorem 10, p.~151].
By this, the sentence
$\underbrace{...~\exists x_i~...}_{\textstyle {x_i \in V}} \Phi$ which is
true in ${\Q}_p$, is also true in ${\Q}_p^{\rm alg}$. Therefore,
for indices $i$ with $x_i \in V$
there exist $w_i \in {\Q}_p^{\rm alg}$
such that ${\Q}_p^{\rm alg} \models \Phi[x_i \mapsto w_i]$.
Since $\Phi$ is quantifier free,
${\Q}_p \models \Phi[x_i \mapsto w_i]$. Thus, by {\bf (9)}
$w_1=r$, so $r \in {\Q}_p^{\rm alg}$.
\vskip 0.2truecm
\par
\noindent
\centerline{\bf 3. Applying R.~M.~Robinson's theorem on definability}
\vskip 0.2truecm
\par
Let a field $\K$ extends $\Q$ and each element
of $\K$ is algebraic over $\Q$. R.~M.~Robinson proved (\cite{Robinson}):
if $r \in \K$ is fixed for all automorphisms of $\K$, then there
exist $U(y),V(y) \in {\Q}[y]$ such that $\{r\}$ is
definable in $\K$ by the formula
$\exists y~(U(y)=0 \wedge x=V(y))$.
By Robinson's theorem $\wt{{\R}^{\rm alg}}={\R}^{\rm alg}$ and
$\wt{{\Q}_p^{\rm alg}}={\Q}_p^{\rm alg}$. Since ${\R}^{\rm alg}$
is an elementary subfield of $\R$ (\cite{Eklof}),
$\wt{\R}=\wt{{\R}^{\rm alg}}$, and
finally $\wt{\R}={\R}^{\rm alg}$. Since ${\Q}_p^{\rm alg}$ is an
elementary subfield of ${\Q}_p$ (\cite{MacIntyre1977},\cite{MacIntyre1986}),
$\wt{{\Q}_p}=\wt{{\Q}_p^{\rm alg}}$,
and finally $\wt{{\Q}_p}={\Q}_p^{\rm alg}$.
\vskip 0.2truecm
\par
\noindent
{\bf Acknowledgement.}  The author thanks the anonymous referee
for valuable suggestions.

Apoloniusz Tyszka\\
Technical Faculty\\
Hugo Ko\l{}\l{}\k{a}taj University\\
Balicka 104, 30-149 Krak\'ow, Poland\\
E-mail address: {\it rttyszka@cyf-kr.edu.pl}
\end{document}